\documentclass[a4paper,12pt]{article}

\usepackage{amsthm}
\usepackage{amsmath,amssymb,latexsym,amsfonts,mathrsfs}

\newtheorem{Thm}{Theorem}[section]
\newtheorem{Prop}{Proposition}[section]

\newcommand{\bN}{\ensuremath{\mathbb{N}}}
\newcommand{\bR}{\ensuremath{\mathbb{R}}}
\newcommand{\bT}{\ensuremath{\mathbb{T}}}
\newcommand{\bZ}{\ensuremath{\mathbb{Z}}}
\newcommand{\cB}{\ensuremath{\mathcal{B}}}
\newcommand{\cF}{\ensuremath{\mathcal{F}}}
\newcommand{\sP}{\ensuremath{\mathscr{P}}}
\newcommand{\sS}{\ensuremath{\mathscr{S}}}

\newcommand{\law}{\stackrel{{\rm law}}{=}}

\newcommand{\eps}{\ensuremath{\varepsilon}}
\newcommand{\e}{{\rm e}} 
\renewcommand{\d}{{\rm d}} 

\newcommand{\tend}[2]{\mathrel{\mathop{\longrightarrow}\limits^{#1}_{#2}}}

\newcommand{\absol}[1]{\left| #1 \right|} 
\newcommand{\rbra}[1]{\left( #1 \right)} 
\newcommand{\cbra}[1]{\left\{ #1 \right\}} 

\newcommand{\proc}[1]{{\it #1}}
\newcommand{\ep}{\qed}
\newcommand{\acks}[1]{{\it Acknowledgements.} #1}

\begin{document}

\begin{center}
{\Large \bf Around Tsirelson's equation, 
or: The evolution process may not explain everything} 
\end{center}
\begin{center}
Kouji \textsc{Yano}\footnote{
Department of Mathematics, Graduate School of Science,
Kobe University, Kobe, JAPAN.}
\quad and \quad 
Marc \textsc{Yor}\footnote{
Laboratoire de Probabilit\'es et Mod\`eles Al\'eatoires, 
Universit\'e Paris VI, Paris, France.}\footnote{
Institut Universitaire de France}
\end{center}

\begin{center}
{\it A survey dedicated to Jean-Paul Thouvenot,} \\
{\it with deep thanks for his constant willingness to listen and help.} 
\end{center}

\begin{abstract}
We present a synthesis of a number of developments 
which have been made around the celebrated Tsirelson's equation (1975), 
conveniently modified in the framework of a Markov chain 
taking values in a compact group $ G $, 
and indexed by negative time. 
To illustrate, we discuss in detail the case of the one-dimensional torus $ G=\bT $. 
\end{abstract}


\section{Introduction}
\ 
\hspace{6pt} 
$ 1^{\circ} $). The contents of this paper, 
which were presented at the Meeting 
``Dynamical Systems and Randomness" 
at IHP, Paris (May 15th, 2009), 
were motivated mainly by the authors' desire 
(\cite{Yor}, \cite{AUY} and \cite{HY}) 
to understand deeply Tsirelson's equation \cite{Cir}. 
This equation shall be discussed in Section \ref{sec: te}, 
while, as a preparation for the main part of the paper, 
we shall discuss the stochastic equation: 
\begin{align}
\eta_k = \xi_k \eta_{k-1} 
\label{eq}
\end{align}
on a compact group $ G $, where $ k $ varies in $ -\bN $, 
$ (\xi_k)_{k \le 0} $ is the ``evolution process", 
and $ (\eta_k)_{k \le 0} $ is the unknown process, 
both taking values in $ G $. 
We believe that equation \eqref{eq} is the ``right" abstraction of Tsirelson's equation, 
as shown in Section \ref{sec: te}. 

More precisely, for every $ k \le 0 $, 
the law of $ \xi_k $ is a given probability $ \mu_k $ on $ G $, 
and the $ \xi_k $'s are assumed independent. 
We shall denote the sequence $ (\mu_k)_{k \le 0} $ simply by $ \mu $. 
It is immediate that 
$ (\eta_k)_{k \le 0} $ is a Markov chain, 
the transitions of which are given by 
\begin{align}
P(\eta_k \in A|\cF^{\eta}_{k-1}) = \mu_k(A \eta_{k-1}^{-1}(\omega)) 
, \quad A \in \cB(G). 
\label{}
\end{align}
The problem is that, as $ k $ varies in $ -\bN $, there is no initial state 
(``at time $ -\infty $"), 
and the study of $ \sP_{\mu} $, the set of the laws of all 
the solutions of \eqref{eq}, necessitates some care. 

We shall also be interested in $ {\rm ex}(\sP_{\mu}) $, 
the set of all extremal points of the compact set $ \sP_{\mu} $, 
as well as in $ \sS_{\mu} $, the set of laws $ P $ 
of ``strong" solutions, i.e.: under $ P $, $ \cF^{\eta}_k \subset \cF^{\xi}_k $, 
hence, since from \eqref{eq}, there is the identity: 
\begin{align}
\xi_k = \eta_k \eta_{k-1}^{-1} 
\label{}
\end{align}
then: $ P \in \sS_{\mu} $ iff: under $ P $, $ \cF^{\eta}_k = \cF^{\xi}_k $, for every $ k $. 

A number of natural questions now arise: 
given $ \mu=(\mu_k)_{k \le 0} $, 
\subitem 
a) is there existence for \eqref{eq}?, 
i.e.: $ \sP_{\mu} \neq \emptyset $. 
\subitem 
b) is there uniqueness?, i.e., $ \sharp (\sP_{\mu}) = 1 $. 
\subitem 
c) is there a strong solution?, i.e.: $ \sS_{\mu} \neq \emptyset $. 

\noindent
We shall see that these different questions may be answered very precisely, 
in particular if $ G = \bT \simeq [0,1) $ is the one-dimensional torus, 
in terms of criteria on $ \mu $.

$ 2^{\circ} $). Interpreting Tsirelson's equation: 
\\ $ \bullet $ 
Before proceeding, we would like to give a ``light" interpretation of equation \eqref{eq}, 
we mean one not to be taken too seriously!: 
$ \eta_k $ describes the ``state of the universe" at time $ k $; 
this state is ``created" by the state at time $ k-1 $, followed 
by the action of the ``evolution" $ \xi_k $. 
The main question is: can today's state of the universe, i.e.: $ \eta_0 $, 
be explained solely from the evolution process?, 
an almost metaphysical question... 
We shall see that the answer(s), 
in terms of $ \mu $, are somewhat paradoxical... 
\\ $ \bullet $ 
Then our ``interpretation" also allows us to justify our quite general choice 
of the probabilities $ (\mu_k)_{k \le 0} $. 
Indeed, today's ``historians of the universe" see the evolution process 
at work, say, for times $ j \in [K,0] $, 
for some large negative $ K $ ($ K $ decreases as ``today" increases...). 
But, they have no knowledge beyond that $ K $, 
hence, we need to make the most general assumptions 
on the $ (\mu_k) $'s to understand all possible cases... 

Thus, depending on the choice of $ \mu $, 
mathematicians give an answer as to 
how much the evolution process determines the present state, 
but this choice remains to be made! 


$ 3^{\circ} $). 
\underline{Plan of the remainder of the paper}: 
\subitem 
\underline{Section \ref{sec: wg}}. A simple example: Wrapping Gaussians on the circle 
\subitem 
\underline{Section \ref{sec: gg}}. The general group framework --- Questions and facts --- 
\subitem 
\underline{Section \ref{sec: te}}. The motivation for this study: Tsirelson's equation 
\subitem 
\underline{Section \ref{sec: fc}}. Some related questions and final comments 

\section{A simple example: Wrapping Gaussians on the circle} \label{sec: wg}

Here, $ G = \bT \simeq [0,1) $; 
$ \eta_k = \exp (2 i \pi \theta_k) $; 
$ \xi_k = \exp (2 i \pi g_k) $, 
with $ \theta_k \in [0,1) $, 
and $ g_k $ Gaussian, centered, $ E[g_k^2] = \sigma_k^2 $. 

The answers to our previous questions are 
radically different depending on whether 
\begin{align}
\sum_{k \le 0} \sigma_k^2 = \infty 
\quad \text{or} \quad 
\sum_{k \le 0} \sigma_k^2 < \infty . 
\label{eq: wrapping gaussian}
\end{align}

\noindent
$ \boxed{\text{The case: $ \textstyle \sum_{k \le 0} \sigma_k^2 = \infty $.}} $ 
\subitem $ \bullet $ 
For any fixed $ k $, $ \eta_k $ is uniformly distributed on the torus, 
(or $ \theta_k $ is uniform on $ [0,1) $), 
independent from the evolution sequence $ (\xi_j)_{j \le 0} $. 
\subitem $ \bullet $ 
$ \sP_{\mu} $ consists of exactly one solution: $ P^*_{\mu} $, 
``the uniform solution". 
\subitem $ \bullet $ 
$ \cF^{\eta}_k = \sigma(\eta_k) \vee \cF^{\xi}_k = \sigma(\eta_j) \vee \cF^{\xi}_k $ 
for $ j \le k $ 
where $ \eta_j $ is independent from $ \cF^{\xi}_k $ and, in fact, even from $ \cF^{\xi}_0 $. 
\subitem $ \bullet $ 
$ \cF^{\eta}_{-\infty } $ is trivial (under $ P^*_{\mu} $).

\noindent
$ \boxed{\text{The case: $ \textstyle \sum_{k \le 0} \sigma_k^2 < \infty $.}} $ 
\subitem $ \bullet $ 
$ \displaystyle \prod_{-N \le j \le 0} \xi_j 
\tend{\rm a.s.}{N \to \infty } \prod_{-\infty }^0 \xi_j $. 
\subitem $ \bullet $ 
Any solution $ (\eta_k) $ satisfies: 
$ \eta_k \tend{\rm a.s.}{k \to -\infty } V $ 
with $ V $ independent of the evolution $ (\xi_j) $, 
\begin{align}
\eta_k = \rbra{ \prod_{-\infty }^k \xi_j } V . 
\label{}
\end{align}
\subitem $ \bullet $ 
$ \cF^{\eta}_k = \cF^{\eta}_{-\infty } \vee \cF^{\xi}_k $ 
and 
$ \cF^{\eta}_{-\infty } = \sigma(V) $. 
\subitem $ \bullet $ 
$ P \in {\rm ex}(\sP_{\mu}) 
\iff P \in \sS_{\mu} 
\iff V=v $, under $ P $ for some constant $ v \in G $. 

Thus, here in the framework of \eqref{eq: wrapping gaussian}, 
there is nonuniqueness iff there is a strong solution, 
which is indeed a puzzling result. 

We give the main arguments of proof when: 
\begin{align}
\sum_{k \le 0} \sigma_k^2 = \infty . 
\label{}
\end{align}
We first show that: 
$ \forall k $, $ \eta_k $ is uniformly distributed, i.e.: 
for $ p \in \bZ $, $ p \neq 0 $, we obtain: 
\begin{align}
\varphi_k(p) 
:=& E[\exp (2 i \pi p \theta_k)] 
\label{} \\
=& E[\exp (2 i \pi p g_k)] \varphi_{k-1}(p) 
\label{} \\
=& \exp \rbra{ - 2 \pi^2 p^2 \rbra{ \sum_{-N \le j \le k} \sigma_j^2 } } \varphi_{-N-1}(p) 
\label{} \\
\to& 0 
\quad (N \to \infty ). 
\label{}
\end{align}
Thus, $ \theta_k $ is uniform on $ [0,1) $, i.e.: $ \eta_k $ is uniform on the torus. 
This reinforces, as we can show, likewise: 
\begin{align}
E[ \exp (2 i \pi p \theta_k) | \xi_k, \xi_{k-1},\ldots,\xi_{-N} ] = 0 . 
\label{}
\end{align}
Then, letting $ N \to \infty $: 
\begin{align}
E[ \exp (2 i \pi p \theta_k) | (\xi_j)_{j \le 0} ] = 0 . 
\label{}
\end{align}
Hence, the independence of $ \eta_k $, for any fixed $ k $, 
from the evolution process $ (\xi_j)_{j \le 0} $. 
In consequence, the law of $ (\eta_k)_{k \le 0} $ is uniquely determined by this 
independence property. 

The triviality of $ \cF^{\eta}_{-\infty } $ will be explained 
(in the next section, Theorem \ref{thm: gen results}) 
by a general result, i.e.: the triviality of $ \cF^{\eta}_{-\infty } $ 
under any $ P \in {\rm ex}(\sP_{\mu}) $, 
but, here, there is only one solution!!

\section{The general group framework --- Questions and facts ---} \label{sec: gg}

Let $ G $ be a general compact group; 
there is the uniform distribution ($ = $ Haar measure), 
and we now take up the discussion of the general questions 
a), b), c) stated in the Introduction. 

Recall that, under $ P \in \sP_{\mu} $, $ (\eta_k)_{k \le 0} $ is a Markov chain, i.e.: 
\begin{align}
P(\eta_k \in A | \cF^{\eta}_{k-1}) = \mu_k (A \eta_{k-1}^{-1}(\omega)) . 
\label{}
\end{align}

\begin{Thm}[Yor \cite{Yor}] \label{thm: unif sol}
For any $ \mu = (\mu_k)_{k \le 0} $, there exists the ``uniform solution" $ P^*_{\mu} $ 
which may be characterized by: $ \forall k $, $ \eta_k $ is uniform, 
independent from the $ (\xi_j)_{j \le 0} $. 
\end{Thm}

\proc{Proof.}
It follows from Kolmogorov's extension theorem, since for 
any $ k $, we may consider the law $ U_k^{(\mu)} $ on $ G^{(-k+1)} $ 
that of $ [\eta_k,\eta_{k+1},\ldots,\eta_0] $; 
then, for $ k<l $ and $ p_{k,l} $, the obvious projection, 
we find that: $ p_{k,l}(U_k^{(\mu)}) = U_l^{(\mu)} $. 
So, the laws $ (U_k^{(\mu)}) $ are consistent, and $ P^*_{\mu} $ exists. 
\ep
\medbreak

Thus, there is always existence; and, we may ask the 2 questions: 
b) is there uniqueness?, i.e.: $ \sharp (\sP_{\mu}) = 1 $?; 
c) does a strong solution exist?, i.e.: $ \sS_{\mu} \neq \emptyset $? 

Let us make the following table: 
\begin{center}
\begin{tabular}{c||c|c}
                 & \multicolumn{2}{c}{$ \boxed{\text{uniqueness}} $} \\
$ \boxed{\text{strong solution}} $ 
                 & holds   & fails \\
\hline
\hline
exists           & $ C_0 $ & $ C_2 $ \\
\hline
does not exist   & $ C_1 $ & $ C_3 $ 
\end{tabular}
\end{center}

\underline{Discussion}: 
We immediately rule out $ C_0 $, 
since, from Theorem \ref{thm: unif sol}, under $ C_0 $, 
the unique solution $ P^*_{\mu} $ is not strong. 
Thus, there remains to discuss the trichotomy: 
$ \boxed{\text{$ C_1 $-$ C_2 $-$ C_3 $}} $ 

We now state several general results: 

\begin{Thm}[Akahori--Uenishi--Yano \cite{AUY}] \label{thm: gen results}
\ 
\subitem {\rm i).} 
$ P (\in \sP_{\mu}) $ is in fact in $ {\rm ex}(\sP_{\mu}) $ 
iff $ \cF^{\eta}_{-\infty } $ is trivial under $ P $. 
\subitem {\rm ii).} 
Let $ g \in G $. The diagonal operator: 
$ \tau_g:(\eta^0_k)_k \to (\eta^0_k g)_k $ acts transitively over $ {\rm ex}(\sP_{\mu}) $, 
i.e.: 
$ \tau_g: P^0 \in {\rm ex}(\sP_{\mu}) \to \tau_g(P^0) \in {\rm ex}(\sP_{\mu}) $ 
and the mapping is surjective. 
\subitem {\rm iii).} 
Any solution $ (\eta_k)_{k \le 0} $ may be represented as: 
\begin{align}
(\eta_k)_k \law (\eta^0_k V)_k 
\label{}
\end{align}
with $ V $ $ G $-valued and independent of $ (\eta^0_k)_k $. 
This yields directly the Krein--Milman integral representation: 
\begin{align}
(\sP_{\mu} \ni) P = \int P(V \in \d v) P^{(\eta^0_k v)_k} 
\label{}
\end{align}
for any given extremal solution $ (\eta^0_k)_k $. 
\end{Thm}

To get a good feeling /introduction/ for the following discussion, 
we recall a result of Csisz\'ar (\cite{Csi}): 
i.e., the ``almost" convergence in law of infinite products of independent random variables. 

\begin{Thm}[Csisz\'ar (\cite{Csi})]
\label{thm: Csiszar}
Let $ (\xi_j)_{j \le 0} $ be our evolution sequence. 
There exists a sequence $ (\alpha _l, l \to -\infty ) $ of deterministic elements of 
$ G $ such that, for any fixed $ k \in -\bN $, the sequence 
\begin{align}
(\xi_k \xi_{k-1} \cdots \xi_l \alpha _l)_{l<k} 
\label{}
\end{align}
converges in law, as $ l \to -\infty $. 
\end{Thm}

We give two illustrations of Theorem \ref{thm: Csiszar}. 

$ 1^{\circ}) $. 
As a first illustration of Theorem \ref{thm: Csiszar}, let us go back to the Gaussian set-up 
of Section \ref{sec: wg}, where we now consider more generally 
\begin{align}
\xi_k = \exp (2 i \pi g_k) , 
\label{}
\end{align}
with $ g_k $ Gaussian, with variance $ \sigma_k^2 $, and mean $ m_k $. 
Then, we may choose 
\begin{align}
\alpha _l = \exp \rbra{ - 2 i \pi \sum_{l \le j \le 0} m_j } 
\label{}
\end{align}
as ``centering sequence". 

$ 2^{\circ}) $. 
To illustrate further Theorem \ref{thm: Csiszar}, or may be, more accurately, 
point 1) of Theorem \ref{thm: HY} below, 
we may consider the case where all 
the laws $ \mu_k $ are the same; then, 
Stromberg \cite{Str} (see also Collins \cite{Col}) 
showed that, for $ \nu $ a given probability on $ G $, 
$ \nu^{*n} $ converges to Haar measure as soon as the smallest subgroup 
which contains the support of $ \nu $ is equal to $ G $.

We may now present a characterization of $ C_1 $ and $ C_2 $. 

\begin{Thm}[Hirayama--Yano \cite{HY}] \label{thm: HY}
The following statements hold: 
\subitem {\rm 1).} 
Uniqueness holds iff, for each $ k \in -\bN $, the products 
$ \xi_k \xi_{k-1} \cdots \xi_l $ converge in law 
as $ l \to -\infty $, to the uniform law on $ G $. 
\subitem {\rm 2).} 
There exists a strong solution iff there exists a sequence 
$ (\alpha _l) $ of deterministic elements of $ G $ such that the products 
\begin{align}
\xi_k \xi_{k-1} \cdots \xi_l \alpha _l 
\quad \text{converge a.s. as $ l \to -\infty $}. 
\label{eq: converge a.s.}
\end{align}
Then, every extremal solution is strong and is the law of the 
a.s. limit of $ (\xi_k \xi_{k-1} \cdots \xi_l \alpha _l g) $, for some $ g \in G $. 
\end{Thm}

Again, to illustrate Theorem \ref{thm: HY}, we may consider 
the general Gaussian hypothesis made after Theorem \ref{thm: Csiszar}; clearly, 
uniqueness holds iff $ \sum_k \sigma_k^2 = \infty $, 
whereas $ C_2 $ holds iff $ \sum_k \sigma_k^2 < \infty $. Note that 
$ C_3 $ never occurs in this set-up.

To give a full discussion of the trichotomy, we go back to 
$ G = \bT \simeq [0,1) $. We introduce: 
\begin{align}
\bZ_{\mu} = \cbra{ p \in \bZ : \text{for some $ k $}, \ 
\prod_{j \le k} \absol{ \int \e^{2 i \pi p x} \mu_j(\d x) } > 0 } . 
\label{}
\end{align}
Then, there is the 

\begin{Prop}[Yor \cite{Yor}]
$ \bZ_{\mu} $ is a subgroup of $ \bZ $; hence, there exists a 
unique integer $ p_{\mu} \ge 0 $ such that $ \bZ_{\mu} = p_{\mu} \bZ $. 
\end{Prop}

This Proposition now allows us to discuss fully 
the trichotomy $ C_1 $-$ C_2 $-$ C_3 $.

\begin{Thm}[Yor \cite{Yor}]
The trichotomy $ C_1 $-$ C_2 $-$ C_3 $ 
may be described as follows, in terms of $ p_{\mu} $: 
\subitem {\rm 1).} 
\underline{Uniqueness in law} iff $ p_{\mu}=0 $, i.e., $ \bZ_{\mu}=\{ 0 \} $. 
Then, 
\subsubitem $ \bullet $ 
$ \cF^{\eta}_{-\infty } $ is trivial; 
\subsubitem $ \bullet $ 
 $ \forall k $, $ \eta_k $ is uniform; 
\subsubitem $ \bullet $ 
$ \cF^{\eta}_k = \sigma(\eta_k) \vee \cF^{\xi}_k $, 
with independence of $ \eta_k $ and $ \cF^{\xi}_k $. 
\subitem {\rm 2).} 
\underline{Existence of a strong solution} iff $ p_{\mu}=1 $, i.e., $ \bZ_{\mu}=\bZ $. 
Then, 
\subsubitem $ \bullet $ 
$ \cF^{\eta}_{-\infty } = \sigma(V) $ 
\subsubitem $ \bullet $ 
for any $ k \in -\bN $, $ \sigma(V) \vee \cF^{\xi}_k $ 
(again, we use the notation of Theorem \ref{thm: gen results}). 
\subitem {\rm 3).} 
\underline{No strong solution, no uniqueness} iff $ p_{\mu} \ge 2 $. 
Then, 
\subsubitem $ \bullet $ 
$ \forall k \in -\bN $, $ \{ p_{\mu} \eta_k \} $ is $ \cF^{\xi}_k $-measurable; 
\subsubitem $ \bullet $ 
$ [ p_{\mu} \eta_k ]/p_{\mu} $ is uniform on 
$ \rbra{ 0,\frac{1}{p_{\mu}},\frac{2}{p_{\mu}},\ldots,\frac{p_{\mu}-1}{p_{\mu}} } $; 
\subsubitem $ \bullet $ 
$ \cF^{\eta}_k = \sigma([ p_{\mu} \eta_k ]) \vee \sigma(\{ p_{\mu} V \}) \vee \cF^{\xi}_k $ 
with independence of the three $ \sigma $-fields. 
\end{Thm}

We may give (many!) sufficient conditions on $ \mu $ which ensure 
either $ C_1 $, or $ C_2 $, or $ C_3 $: 
\subitem {\rm a).} 
There is uniqueness in law (i.e.: $ p_{\mu}=0 $) as soon as 
$ \xi_{n_j} \law \exp ( 2 \pi i \eps_j \gamma) $ 
for some subsequence $ (n_j) $, 
some $ \eps_j \in \bR $ with $ |\eps_j| \tend{}{(j \to \infty )} \infty $, 
and some random variable $ \gamma $ with absolutely continuous density. 
\subitem {\rm b).} 
Assume $ \mu_j = \nu $, for all $ j $. Then: 
\subsubitem {\rm i).} 
$ p_{\mu}=0 $ iff $ \nu $ is not arithmetic; 
\subsubitem {\rm ii).} 
$ p_{\mu}=1 $ iff $ \exists x \in \bR $, $ \nu(x+\bZ)=1 $; 
\subsubitem {\rm iii).} 
$ p_{\mu}=p \ge 2 $ if $ \nu $ charges precisely 
$ \rbra{ 0,\frac{1}{p},\frac{2}{p},\ldots,\frac{p-1}{p} } $. 

\proc{Remark.}
A special case of a) above is when 
$ \eps_j = j $ 
for some random variable $ \gamma $ with absolutely continuous density. 
This is called ``Poincar\'e roulette wheel, leading to equidistribution"; 
see \cite[Theorem 3.2]{Eng} for detail. 
\medbreak

\section{The motivation for this study: Tsirelson's equation} \label{sec: te}

A result of Zvonkin \cite{Zvo} (see also Zvonkin--Krylov \cite{ZK}) 
asserts that the stochastic differential equation 
driven by BM: 
\begin{align}
X_t = B_t + \int_0^t \d s \ b(X_s) , 
\label{}
\end{align}
where $ b(\cdot) $ is only assumed to be bounded and Borel 
enjoys {\em strong uniqueness}. 
(For $ h $ such that $ \frac{1}{2} h'' + b h' = 0 $, 
the process $ Y_t = h(X_t) $ solves 
$ Y_t = \int_0^t (h' \circ h^{-1})(Y_s) \d B_s $ 
where $ h' \circ h^{-1} $ is Lipschitz.) 

Then, the question arose whether the same strong uniqueness result 
might still be true with a bounded Borel drift 
depending more generally on the past of $ X $, i.e.: 
\begin{align}
X_t = B_t + \int_0^t \d s \ b(X_u,u \le s) 
\label{eq: 22}
\end{align}
(Uniqueness in law is ensured by Girsanov's theorem). 
Tsirelson gave a negative answer to this question by producing the drift: 
\begin{align}
b(X_u,u \le s) = \sum_{k \in -\bN} \cbra{ \frac{X_{t_k}-X_{t_{k-1}}}{t_k - t_{k-1}} } 
1_{(t_k,t_{k+1}]}(s) 
\label{eq: 23}
\end{align}
for any sequence $ t_k \downarrow 0 $ as $ k \downarrow -\infty $. 

To prove that the solution is non-strong, it suffices to study the 
discrete time skeleton equation: 
\begin{align}
\frac{X_{t_{k+1}}-X_{t_{k}}}{t_{k+1} - t_{k}} 
= 
\frac{B_{t_{k+1}}-B_{t_{k}}}{t_{k+1} - t_{k}} 
+ 
\cbra{ \frac{X_{t_k}-X_{t_{k-1}}}{t_k - t_{k-1}} } 
\label{}
\end{align}
i.e.: 
\begin{align}
\eta_k = \xi_k + \{ \eta_{k-1} \} . 
\label{}
\end{align}
Slight modifications of our previous arguments show that: 
$ \forall k $, $ \{ \eta_k \} $ is independent from the BM, and uniformly 
distributed on $ [0,1) $; moreover, $ \forall t $, $ \forall t_k \le t $, 
\begin{align}
\cF^X_t = \sigma(\{ \eta_{k-1} \}) \vee \cF^{B}_t . 
\label{eq: 26}
\end{align}

For many further references, see Tsirelson's web page \cite{Tsi}.

\section{Some related questions and final comments} \label{sec: fc} 
\ 
\hspace{6pt} 
a). The case $ C_1 $ [only $ P^*_{\mu} $ solution] 
gives a beautiful example where: 
\begin{align}
\cF^{\eta}_k 
=& \cF^{\eta}_j \vee \cF^{\xi}_k 
, \quad j \le k 
\label{} \\
=& \bigcap_j (\cF^{\eta}_j \vee \cF^{\xi}_k) 
\neq \rbra{ \bigcap_j \cF^{\eta}_j } \vee \cF^{\xi}_k 
\label{}
\end{align}
since in the $ C_1 $ case, $ \bigcap_j \cF^{\eta}_j = \cF^{\eta}_{-\infty } $ is trivial. 
This is a discussion which has been a trap 
for a number of very distinguished mathematicians... 
See, e.g., N.~Wiener (\cite{Wie}, Chap 2) 
where he assumes that $ \cap $ and $ \vee $ may be interverted for 
$ \sigma $-fields; if so, from Wiener's set-up, this would lead to: 
$ K $-automorphism is always Bernoulli! 

b). 
As \eqref{eq: 26} above clearly shows, 
any solution to Tsirelson's equation \eqref{eq: 22}-\eqref{eq: 23} 
is not strong, i.e., $ (X_t) $ cannot be recovered from the Brownian motion 
$ (B_t) $ in \eqref{eq: 22}-\eqref{eq: 23}. 
Nevertheless, as shown in Emery--Schachermayer \cite{ES}, 
the natural filtration of $ X $, that is: 
$ (\cF^X_t)_{t \ge 0} $ is generated by \underline{some} Brownian motion $ (\beta_t)_{t \ge 0} $; 
thus, in that sense, the filtration $ (\cF^X_t)_{t \ge 0} $ is a strong Brownian filtration. 
The question then arose naturally whether under any probability $ Q $ on $ C([0,1],\bR) $, 
equivalent to Wiener measure, the natural filtration of the canonical coordinate process 
is always a strong Brownian filtration. 
This is definitely not the case, as shown, e.g., 
in Dubins--Feldman--Smorodinsky--Tsirelson \cite{DFST}. 
The paper \cite{ES} contains a number of important references 
around the topic treated in \cite{DFST}, 
where again the role of B.~Tsirelson has been crucial. 
In particular, it is another beautiful result of B.~Tsirelson \cite{TsiGAFA} 
that the natural filtration of the Brownian spider with $ N (\ge 3) $ legs 
is not strongly Brownian.

c). 
Related to the problem b), we may ask the following question: 
For a solution $ (\eta_k)_{k \le 0} $ of equation \eqref{eq}, 
when does there exist \underline{some} sequence $ \theta = (\theta_k)_{k \le 0} $ 
of independent random variables such that 
$ \eta_k \in \cF^{\theta}_k $ for any $ k \le 0 $? 
There are lots of studies in the case where $ \eta $ is a stationary process; 
some positive answers to this question are found in \cite{Ros60} and \cite{Hanson1}. 
See \cite{Lau} for historical remarks and related references.

d). 
For stationary processes with common law $ \mu $ 
on a state space $ S $ with a continuous group action of a locally compact group $ G $, 
Furstenberg \cite[Definition 8.1]{F} introduced the notion of a ``$ \mu $-boundary". 
When we confine ourselves to the case where $ S=G $ with the canonical group action 
and where the noise process is assumed to be identically distributed, 
Furstenberg's $ \mu $-boundary 
is essentially the same as a strong solution in our terminology, 
and Theorem 14.1 of \cite{F} coincides with our Theorem \ref{thm: HY}. 

Let $ G $ be a compact group with countable basis 
endowed with metric $ d $ and let $ S $ be a compact space with continuous $ G $-action. 
Furstenberg (\cite{F2}), in his study of stationary measures, 
proved the following: 
If the action is {\em distal}, i.e., 
for any $ x,y \in S $ with $ x \neq y $, it holds that $ \inf_{g \in G} d(gx,gy) > 0 $, 
then the action is {\em stiff}, i.e., 
for any probability law $ \mu $ on $ G $ whose support generates $ G $, 
any $ \mu $-invariant probability measure on $ S $ is $ G $-invariant. 
Coming back to our setting, we assume $ S=G $ with the canonical group action. 
Then the well-known theorem by Birkhoff \cite{B} and Kakutani \cite{K} 
shows that we may choose the metric $ d $ so that it remains invariant under the action of $ G $, 
i.e., $ d(gx,gy) = d(x,y) $ for any $ g \in G $. 
Thus the canonical action of $ G $ may be assumed to be distal, 
and we see that, then, the canonical action is stiff.

e). 
We would also like to point out the relevance of the 
It\^o--Nisio \cite{IN} study of all stationary solutions $ (X_t)_{-\infty <t<\infty } $ 
of some stochastic differential equations 
driven by Brownian motion $ (B_t)_{-\infty <t<\infty } $. 
They discuss whether either of the following properties holds: 
{\subitem (i) 
$ \cF^X_{-\infty ,t} \subset \cF^B_{-\infty ,t} $; 
\subitem (ii) 
$ \cap_{t} \cF^X_{-\infty ,t} $ is trivial; 
\subitem (iii) 
$ \cF^X_{-\infty ,t} \subset \cF^X_{-\infty ,s} \vee \cF^B_{s,t} $ for $ s<t $,} 

\noindent
where we have used obvious notations for $ \sigma $-fields. 
It\^o--Nisio \cite[Section 13]{IN} identify cases where (ii) and (iii) hold, 
but not (i); 
this is similar to the situation in Tsirelson's original equation \eqref{eq: 22}-\eqref{eq: 23}, 
or more generally the $ C_1 $ case. 
But, even worse, It\^o--Nisio \cite[Section 14]{IN} also discuss a case 
where (iii) does not hold. This case originates from Girsanov's equation \cite{G} 
with diffusion coefficient $ |x|^{\alpha } $ for $ \alpha <1/2 $, 
which is well-known to generate non-uniqueness in law.

\begin{center}
\rule{50mm}{0.5pt}
\end{center}

As a word of conclusion, 
although we do not claim that equation \eqref{eq} has a deep ``cosmological value", 
it would probably never come to the mind of ``universe historians" that, 
in some cases, despite ``the emptiness of the beginning", 
today's state may be independent of the evolution mechanism. 
Thus, seen in this light, Tsirelson's equation, and its abstraction \eqref{eq} 
provide us with a beautiful and mind-boggling statement. 

\acks{We are grateful to Professor Benjamin Weiss for making several deep comments 
on a first draft, 
and to Professor Bernard Roynette 
for providing us with the references \cite{Col} and \cite{Str}.}

\def\cprime{$'$} \def\cprime{$'$}

\end{document}